\documentclass[12pt]{article}
\usepackage{amsfonts}

\usepackage{amssymb,amsmath}
\date{ }
\def\ma{\mathbb}
\def\BE {\begin{eqnarray}}
\def\EE {\end{eqnarray}}
\def\BC {\begin{eqnarray*}}
\def\EC {\end{eqnarray*}}
\def\OPLUS#1{\raisebox{-5pt}{\mbox{$\begin{array}{c}
\oplus\\[-5pt]\scriptstyle#1
\end{array}$}}}

\title{{\bf\Large The Derivation algebra and automorphism group of the generalized Ramond $N=2$ superconformal algebra}
\thanks{This work is supported by NSF grants 10726057 of China.}}

\author{
{{\bf Jiayuan Fu} ~~ {\bf Yongcun Gao} }
\\ {\small {\it Department of Mathematics,  Communication University of China,}}
 \\
  {\small {\it Beijing 100024, China }}
  \\{\small {\it Email: fujy@cuc.edu.cn}}
}

\begin{document}
\maketitle
\begin{center}{\large \bf Abstract}\end{center}
  \hskip 20pt{\small
   In this paper, we give the definition of the generalized
   Ramond $N=2$ superconformal algebras and discuss the derivation
   algebra and the automorphism group. }
 \vskip 0.6 cm
\begin{center}{\large \bf 1. Introduction}\end{center}
\def\theequation{1.\arabic{equation}}

Superconformal algebras have been constructed for almost three
decades(\cite {A1, K1}). Since then the study of superconformal
algebras has made much progress in both mathematics and physics. Kac
and van de Leuer (\cite{KL}), Cheng and Kac (\cite{CK})  have
classified all possible superconformal algebras and then Kac has
proved that their classification is complete.

Structure theory and representation theory are also the most
important two theories on Lie superalgebras. For $N=2$
superconformal algebras, remarkable efforts have been taken by
several research groups (\cite{BF, D1, D2, DK, BR, KE}). In
\cite{FJS}, the modules of intermediate series over Ramond $N=2$
superconformal algebra has been classified. For the Ramond $N=2$
superconformal algebra, its even part is the so called twisted
Heisenberg-Virasoro algebra (\cite{SJ}). The derivation algebras of
infinite-dimensional Heisenberg algebra, Virasoro algebra and the
twisted Heisenberg-Virasoro algebra are given by Jiang(\cite{J}),
Zhu, Meng(\cite{ZM}) and Shen, Jiang (\cite{SJ}) respectively. Su,
Song and Zhou have given the structures of derivation algebra of
Weyl and  Block type algebras in {\cite{SS, S1, S2, SZ}}. Most of
those algebras are finitely generated. In this paper, we will give
the derivation algebra and the automorphism group of the generalized
Ramond $N=2$ superconformal algebra. In general, this algebra  is
infinitely generated.

Let ${\ma F}$ be a field of characteristic zero, $\Gamma$ be an
additive subgroup of ${\ma F}$, ${\bf 0}$ be the identity element of
$\Gamma$. The generalized Ramond $N=2$ algebra is defined as below:
$${\cal L}={\cal L}_{\overline{0}}\oplus{\cal L}_{\overline{1}},$$
where \BC {\cal L}_{\overline{0}}=\mbox{span}_{\ma C}\{ L_\alpha,
H_\beta, c\mid \alpha, \beta\in\Gamma \},  ~~~  {\cal
L}_{\overline{1}}=\mbox{span}_{\ma C}\{G_\alpha^-, G_\beta^+\mid
\alpha, \beta\in\Gamma \},\EC with the following relations:
\BE\label{a23}
\begin{array}{lll}
 [ L_\alpha, L_\beta ] & = &
 (\alpha-\beta)L_{\alpha+\beta}+\frac{c}{12}(\alpha^3-\alpha)\delta_{\alpha+\beta,0}, \\ [0.3 cm]
   [ H_\alpha, H_\beta ] &=& \frac{\alpha}{3} c \delta_{\alpha+\beta,0}, \\ [0.3 cm]
    [ L_\alpha, H_\beta ] &=& -\beta H_{\alpha+\beta}, \\[0.3 cm] [ L_\alpha, G_\beta^{\pm} ]&
  = & (\frac{\alpha}{2}-\beta)G_{\alpha+\beta}^{\pm}, \\[0.3 cm]
   [ H_\alpha, G_\beta^{\pm} ] &=&  \pm G_{\alpha+\beta}^{\pm}, \\ [0.3 cm]
    [ G_\alpha^-, G_\beta^+ ]&=&
  2L_{\alpha+\beta}-(\alpha-\beta)H_{\alpha+\beta}+\frac{c}{3}(\alpha^2-\frac{1}{4})\delta_{\alpha+\beta,0},\\
 \end{array} \EE
where $\delta_{m,n}$ is the Kronecker notation which satisfies that $\delta_{m,n}=\left\{ \begin{array}{ll} 1, &  m=n \\
0, &  m\not= n \end{array}\right.$.
 Throughout this paper, ${\cal L}$ denotes the {\bf generalized Ramond $N=2$ superconformal
 algebra}. It is easy to see that ${\cal L}$ is a $\Gamma$-graded
 algebra:
 $${\cal L}=\OPLUS{\alpha\in \Gamma}{\cal L}_\alpha,$$
 where ${\cal L}_\alpha=\{x\in{\cal L}\mid [L_0, x]=-\alpha
 x\}=$span$_{\ma F}\{L_\alpha, H_\alpha, G^\pm_\alpha, \delta_{\alpha, 0}c\}$.

In the following section, we will discuss the derivation algebra of
${\cal L}$. Our main result in this section is theorem 2.3. Then the
automorphism group will be given in the last section, and we will
obtain the main theorem 3.4.  \vskip 0.8 cm
\begin{center}{\bf 2. The derivation algebra of ${\cal L}$}\end{center}
\def\theequation{2.\arabic{equation}}
\setcounter{equation}{0}
 \mbox{} \hskip 10pt

In this section, we denote the derivation algebra of ${\cal L}$ by
$Der{\cal L}$, the set of inner derivations by $ad{\cal L}$. Since
${\cal L}$ is a ${\ma Z}_2$-graded algebra, then $$ Der{\cal
L}=(Der{\cal L})_{\overline{0}}\oplus(Der{\cal L})_{\overline{1}},$$
where
$$(Der{\cal L})_{\overline{0}}=\{D\in Der{\cal L}\mid D({\cal
L}_{\overline{i}})\subseteq{\cal L}_{\overline{i}}, ~
\overline{i}\in {\ma Z}_2 \},$$  denotes the set of even derivations
of ${\cal L}$, and $$(Der{\cal L})_{\overline{1}}=\{D\in Der{\cal
L}\mid D({\cal L}_{\overline{i}})\subseteq{\cal L}_{\overline{i+1}},
~ \overline{i}\in {\ma Z}_2 \},$$ denotes the set of odd derivations
of ${\cal L}$. Note that ${\cal L}$ is $\Gamma$-graded,
  we have that $Der{\cal L}$ is also $\Gamma$-graded:
define $$ (Der{\cal L})_\gamma=\{D\in Der{\cal L}\mid D({\cal
L}_\beta)\subseteq {\cal L}_{\beta+\gamma}, \beta\in\Gamma\}.
$$
Then we have that $Der{\cal L}=\OPLUS{\gamma\in\Gamma}(Der{\cal
L})_\gamma$ (see {\cite{S1,S2,SS,SZ}}). Obviously, the
$\Gamma$-graded and ${\ma Z}_2$-graded structures are compatible,
i.e.,  $Der{\cal L}_{\overline{0}}=\OPLUS{\gamma\in\Gamma}(Der{\cal
L})'_\gamma$, $Der{\cal
L}_{\overline{1}}=\OPLUS{\gamma\in\Gamma}(Der{\cal L})''_\gamma$,
where $(Der{\cal L})'_\gamma\oplus(Der{\cal L})''_\gamma=(Der{\cal
L})_\gamma$,  the ${\ma Z}_2$-graded structure of $(Der{\cal
L})_\gamma$ is similar. In the following, we will only discuss the
homogenous derivation $D$. \vskip 0.3cm {\bf Lemma 2.1.}\quad The
odd derivations of ${\cal L}$ are all inner derivations, i.e.,
$(Der{\cal L})_{\overline{1}}=ad {\cal L}_{\overline{1}}$. \vskip
0.2 cm
 {\bf Proof.}\quad Obviously, we
only need to prove that $ (Der{\cal L})_{\overline{1}}\subseteq ad
{\cal L}_{\overline{1}}$.  If $D\in (Der{\cal
L})_{\overline{1}}\cap(Der{\cal L})_\gamma$, $\gamma\in\Gamma$,
suppose that \BC\begin{array}{lll} D(L_\alpha) &=&
a_{\gamma,\alpha}G^+_{\alpha+\gamma}+b_{\gamma,\alpha}G^-_{\alpha+\gamma},
\\ [0.2cm] D(H_\alpha)&=& c_{\gamma,\alpha}G^+_{\alpha+\gamma}+d_{\gamma,\alpha}G^-_{\alpha+\gamma},
\\ [0.2cm]
D(G^+_\alpha)&=&
e_{\gamma,\alpha}L_{\alpha+\gamma}+f_{\gamma,\alpha}H_{\alpha+\gamma}+m_\alpha\delta_{\alpha+\gamma,0}c,\\
[0.2cm]D(G^-_\alpha)&=&
g_{\gamma,\alpha}L_{\alpha+\gamma}+h_{\gamma,\alpha}H_{\alpha+\gamma}+n_\alpha\delta_{\alpha+\gamma,0}c.\end{array}\EC
Since $D$ is a derivation, by the following equations:
\BE\begin{array}{lll} \label{g2} D[H_\alpha, H_\beta] & = &
[D(H_\alpha), H_\beta]+[H_\alpha, D(H_\beta)],\\ [0.3 cm]
\label{g3}D[L_\alpha, H_\beta] & = & [D(L_\alpha),
H_\beta]+[L_\alpha, D(H_\beta)],\\ [0.3 cm]  \label{g4}D[H_\alpha,
G^\pm_\beta]& = & [D(H_\alpha), G^\pm_\beta]+[H_\alpha,
D(G^\pm_\beta)],\end{array}\EE for all $\alpha, \beta\in\Gamma$, we
can obtain some relations between these coefficients:
 $$ \begin{array}{rll}
c_{\gamma,\alpha}& = & c_{\gamma,\beta}, \\ [0.2cm]
d_{\gamma,\alpha} & = & d_{\gamma,\beta},  \\ [0.2cm]
e_{\gamma,\alpha+\beta}& = & 2d_{\gamma,\alpha},  \\ [0.2 cm]
f_{\gamma,\alpha+\beta} & = &
\alpha e_{\gamma,\beta}-(\alpha+\gamma-\beta)d_{\gamma,\alpha}, \\
[0.2 cm] m_{\alpha+\beta} & = &
\frac{1}{3}(\beta^2-\frac{1}{4})d_{\gamma,\alpha}
+\frac{\alpha}{3}f_{\gamma,\beta}, ~~~ \alpha+\beta+\gamma=0, \\
[0.2 cm] g_{\gamma,\alpha+\beta} & = & -2c_{\gamma,\alpha}, \\ [0.2
cm] h_{\gamma,\alpha+\beta} & = & -\alpha
g_{\gamma,\beta}-(\alpha+\gamma-\beta)c_{\gamma,\alpha},\\ [0.2 cm]
n_{\alpha+\beta} & = &
-\frac{1}{3}(\beta^2-\frac{1}{4})c_{\gamma,\alpha}
-\frac{\alpha}{3}h_{\gamma,\beta}, ~~~ \alpha+\beta+\gamma=0, \\
[0.2 cm]\beta c_{\gamma,\alpha+\beta} & = &
a_{\gamma,\alpha}-(\frac{\alpha}{2}-\beta-\gamma)c_{\gamma,\beta}, \\
[0.2 cm] \beta d_{\gamma,\alpha+\beta} & = &
-b_{\gamma,\alpha}-(\frac{\alpha}{2}-\beta-\gamma)d_{\gamma,\beta}.
\end{array}$$
It is not difficult to see that $c_{\gamma,\alpha}$ and
$d_{\gamma,\alpha}$ are contents. Set $c_{\gamma,\alpha}=\xi_0,
d_{\gamma,\alpha}=\xi_1$, for all $\alpha\in\Gamma$, where $\xi_0,
\xi_1\in{\ma F}$, then we can obtain the following results:
$$\begin{array}{rllrll} a_{\gamma,\alpha} &= &
(\frac{\alpha}{2}-\gamma)\xi_0, & b_{\gamma,\alpha} & = &
-(\frac{\alpha}{2}-\gamma)\xi_1, \\ [0.2 cm] e_{\gamma,\alpha} &= &
2\xi_1, & f_{\gamma,\alpha} & = & (\alpha-\gamma)\xi_1, \\
[0.2 cm] g_{\gamma,\alpha} &= &
-2\xi_0, & h_{\gamma,\alpha} & = & (\alpha-\gamma)\xi_0, \\
[0.2 cm] m_\alpha & = & \frac{1}{3}(\alpha^2-\frac{1}{4})\xi_1, &
n_\alpha & = & -\frac{1}{3}(\alpha^2-\frac{1}{4})\xi_0.
\end{array}$$
By these equations, we can easily check that
$D=ad(-\xi_0G^+_\gamma+\xi_1G^-_\gamma)$. \hfill $\Box$ \vskip 0.3cm
{\bf Lemma 2.2.}\quad If $D\in (Der{\cal
L})_{\overline{0}}\cap(Der{\cal L})_\gamma$, and $\gamma\not={\bf
0}$, then $D\in ad{\cal L}$. \vskip 0.2 cm
 {\bf Proof.}\quad
$D(c)=0$, since $\gamma\not={\bf 0}$. Without confusion, we use the
same symbols: suppose that \BC\begin{array}{lll} D(L_\alpha) &=&
a_{\gamma,\alpha}L_{\alpha+\gamma}+b_{\gamma,\alpha}H_{\alpha+\gamma}+m_\alpha\delta_{\alpha+\gamma,0}c,
\\ [0.2cm] D(H_\alpha)&=& c_{\gamma,\alpha}L_{\alpha+\gamma}+d_{\gamma,\alpha}H_{\alpha+\gamma}+n_\alpha\delta_{\alpha+\gamma,0}c,
\\ [0.2cm]
D(G^+_\alpha)&=&
e_{\gamma,\alpha}G^+_{\alpha+\gamma}+f_{\gamma,\alpha}G^-_{\alpha+\gamma},\\
[0.2cm]D(G^-_\alpha)&=&
g_{\gamma,\alpha}G^+_{\alpha+\gamma}+h_{\gamma,\alpha}G^-_{\alpha+\gamma}.\end{array}\EC
Also by (\ref{g2}) and \BE\label{g6} D[L_\alpha, G^\pm_\beta]=
[D(L_\alpha), G^\pm_\beta]+[L_\alpha, D(G^\pm_\beta)],\EE
\BE\label{g7} D[G^-_\alpha, G^+_\beta]= [D(G^-_\alpha),
G^+_\beta]+[G^-_\alpha, D(G^+_\beta)],\EE we have that
 \BC  \label{g10}\beta c_{\gamma,\alpha}& = & \alpha c_{\gamma,\beta}, \\
\label{g14} e_{\gamma,\alpha+\beta}& = & (\frac{\alpha+\gamma}{2}-\beta)c_{\gamma,\alpha}+d_{\gamma,\alpha}+e_{\gamma,\beta},  \\
 \label{g16}  f_{\gamma,\alpha+\beta} & = & -f_{\gamma,\beta},
 \hskip 1.5 cm
   g_{\gamma,\alpha+\beta} ~ = ~ -g_{\gamma,\beta},\\
 \label{g20} h_{\gamma,\alpha+\beta}& = & -(\frac{\alpha+\gamma}{2}-\beta)c_{\gamma,\alpha}+d_{\gamma,\alpha}+h_{\gamma,\beta},  \\
 \label{g22} -\beta c_{\gamma,\alpha+\beta} & = &
 (\alpha-\beta-\gamma)c_{\gamma,\beta} \\
 \label{g26} (\frac{\alpha}{2}-\beta)e_{\gamma,\alpha+\beta} &=& (\frac{\alpha+\gamma}{2}-\beta)
a_{\gamma,\alpha}+b_{\gamma,\alpha}
+(\frac{\alpha}{2}-\beta-\gamma)e_{\gamma,\beta},\\
\label{g28}  (\frac{\alpha}{2}-\beta)h_{\gamma,\alpha+\beta} &=&
(\frac{\alpha+\gamma}{2}-\beta)a_{\gamma,\alpha}-b_{\gamma,\alpha}
+(\frac{\alpha}{2}-\beta-\gamma)h_{\gamma,\beta},\\ \label{g30}
-\beta n_{\alpha+\beta}&=&
\frac{\alpha+\gamma}{3}b_{\gamma,\alpha}+\frac{\alpha^3-\alpha}{12}c_{\gamma,\beta},
~~ \alpha+\beta+\gamma=0, \\
 \label{g32} 2a_{\gamma,\alpha+\beta}-(\alpha-\beta)c_{\gamma,\alpha+\beta} & = & 2h_{\gamma,\alpha}+2e_{\gamma,\beta},\\
\label{g34}
2b_{\gamma,\alpha+\beta}-(\alpha-\beta)d_{\gamma,\alpha+\beta} & = &
-(\alpha+\gamma-\beta)h_{\gamma,\alpha}-(\alpha-\beta-\gamma)e_{\gamma,\beta}, \\
\label{g36} 6m_{\alpha+\beta}-3(\alpha-\beta)n_{\alpha+\beta} & = &
(\beta^2-\frac{1}{4})h_{\gamma,\alpha}+(\alpha^2-\frac{1}{4})e_{\gamma,\beta},
  \alpha+\beta=-\gamma. \EC

From these equations, we can obtain that for all $\alpha\in\Gamma$
$$ c_{\gamma,\alpha}=f_{\gamma,\alpha}=g_{\gamma,\alpha}=0, ~
b_{\gamma,\alpha}=\eta_\gamma, ~
n_{\alpha}=\frac{1}{3}\eta_\gamma,$$ where $\eta_\gamma\in{\ma F},$
and \BC\begin{array}{llllll} a_{\gamma,\alpha} & = &
\gamma^{-1}(\gamma-\alpha)(2h_{\gamma,0}+2\gamma^{-1}\eta_\gamma), &
d_{\gamma,\alpha} & = &
-\gamma^{-1}\alpha(2h_{\gamma,0}+2\gamma^{-1}\eta_\gamma), \\ [0.2
cm] e_{\gamma,\alpha} & = &
\gamma^{-1}(\gamma-2\alpha)h_{\gamma,0}+2\gamma^{-2}(\gamma-\alpha)\eta_\gamma,
& h_{\gamma,\alpha} & = &
\gamma^{-1}(\gamma-2\alpha)h_{\gamma,0}-2\gamma^{-2}\alpha\eta_\gamma,
\\ [0.2 cm]
m_\alpha & = &
\frac{1}{6}\gamma^{-1}(\alpha-\alpha^3)(h_{\gamma,0}+\gamma^{-1}\eta_\gamma).
& & & \end{array}\EC It is not difficult to see that, if $D\in
(Der{\cal L})_{\overline{0}}\cap(Der{\cal L})_\gamma$, and
$\gamma\not=0$,  $$ D=ad(k_1L_\gamma+k_2H_\gamma),
$$ where $k_1=2\gamma^{-1}(h_0+\gamma^{-1}\eta_\gamma),
k_2=\gamma^{-1}\eta_\gamma$, $h_0,\eta_\gamma\in{\ma F}$.
\hfill$\Box$

If $D\in (Der{\cal L})_{\overline{0}}\cap(Der{\cal L})_0$, the case
is a little difference, and we should note that $D(c)=i_0c$,
$i_0\in{}\ma F$. By (\ref{g2})-(\ref{g7}), used the similar method,
we can obtain that for all $\alpha, \beta\in\Gamma$. $$
n_{\gamma,\alpha}=
b_{\gamma,\alpha}=c_{\gamma,\alpha}=f_{\gamma,\alpha}=g_{\gamma,\alpha}=i_0=m_0=0,$$
 and
$$ \left \{ \begin{array}{ll}
h_{\gamma,\alpha}=e_\alpha-2e_0, ~~
a_{\gamma,\alpha}=d_{\gamma,\alpha}=e_\alpha-e_0,
\\ [0.2 cm] a_{\gamma,\alpha+\beta}=a_{\gamma,\alpha}+a_{\gamma,\beta}. \end{array}\right.$$
 Denote by Hom$_{\ma Z}(\Gamma,
{\ma F})$ the set of additive group homomorphisms from $\Gamma$ to
${\ma F}$.
 For
$\varphi \in$Hom$_{\ma Z}(\Gamma, {\ma F})$, we define a derivation
also denoted $\varphi$ by $\varphi
(x_\alpha)=\varphi(\alpha)x_\alpha$, for all $x_\alpha\in{\cal L}$.
We have the following theorem: \vskip 0.3cm {\bf Theorem 2.3.}\quad
$Der{\cal L}$ is spanned by $ad{\cal L}$ and Hom$_{\ma Z}(\Gamma,
{\ma F})$. In details, $Der{\cal L}=\OPLUS{\alpha\in\Gamma}(Der{\cal
L})_\alpha =(Der{\cal L})_{\overline{0}}\oplus(Der{\cal
L})_{\overline{1}},$ where
$$(Der{\cal L})_{\overline{i}}=\left\{\begin{array}{ll} ad{\cal
L}_{\overline{1}}, & \overline{i}=\overline{1}, \\ [0.2 cm] ad{\cal
L}_\alpha, & \overline{i}=\overline{0}, \alpha\not=0,
\\ [0.2 cm] ad{\cal L}'_0+\mbox{Hom}_{\ma Z}(\Gamma, {\ma F}), &
\overline{i}=\overline{0}, \alpha=0, \end{array}\right.$$ where
${\cal L}'_0=$span$_{\ma F}\{L_0, H_0\}$, and $Der{\cal
L}_{\overline{0}}\cap ad{\cal L}_0=\{ad(kL_0+lH_0)\mid k, l\in{\ma
F}\}$. \vskip 0.2 cm {\bf Proof.}\quad By the above argument, we
know that if $D\in (Der{\cal L})_{\overline{0}}\cap(Der{\cal L})_0$,
then \BC\begin{array}{lll} D(L_\alpha)=\varphi(\alpha) L_\alpha, &
D(H_\alpha)=\varphi(\alpha) H_\alpha, &
D(G^+_\alpha)=(\varphi(\alpha)+e_0)G^+_\alpha, \\ [0.2 cm]
D(G_\alpha^-)=(\varphi(\alpha)-e_0)G^-_\alpha, & D(c)=0,
\end{array}\EC where $\varphi\in$Hom$_{\ma Z}(\Gamma, {\ma F})$,
i.e. $\varphi(\alpha+\beta)=\varphi(\alpha)+\varphi(\beta)$ for any
$\alpha, \beta\in\Gamma$, and $e_0\in{\ma F}$. Obviously, $D\in
ad{\cal L}$ if and only if $\varphi(\alpha)=k\alpha$, $k\in{\ma F}$.
And at this time, $D=ad(-kL_0+e_0H_0)$.
 \hfill$\Box$
 \vskip 0.6 cm
\begin{center}{\bf 3. The automorphism group of ${\cal L}$}
\end{center}
\def\theequation{3.\arabic{equation}}
\setcounter{equation}{0} \mbox{}

Denote by $Aut{\cal L}$ the automorphism group. By the relations of
${\cal L}$, we have that $G^+_\alpha, G^-_{\beta}$ are locally
ad-nilpotent elements for any $\alpha, \beta\in\Gamma$. In this
section, we denote ${\ma F}^*$ the multiply group which is generated
by the non-zero elements of ${\ma F}$.  \vskip 0.2 cm {\bf Lemma
3.1.}\quad For any $\sigma\in Aut{\cal L}$, $\sigma(G^\pm_\gamma)\in
$ span$_{\ma F}\{G^+_\alpha \mid \alpha\in\Gamma\}\cup$ span$_{\ma
F}\{ G^-_\alpha\mid \alpha\in\Gamma \}$. \vskip 0.2 cm
 {\bf
Proof.}\quad  Note that for any $\sigma\in Aut{\cal L}$,
$\sigma({\cal L}_{\overline{i}})={\cal L}_{\overline{i}}$,
$\overline{i}\in{\ma Z}_2$, and $\sigma(x)$ is also a locally
ad-nilpotent element if $x$ is a locally ad-nilpotent element, and
$G^+_\alpha$ and $G^-_\beta$ are not abelian, then it is not
difficult to prove that $\sigma(G^\pm_\gamma)\in $ span$_{\ma
F}\{G^+_\alpha \mid \alpha\in\Gamma\}\cup$ span$_{\ma F}\{
G^-_\alpha\mid \alpha\in\Gamma \}$. \hfill$\Box$ \vskip 0.3 cm {\bf
Lemma 3.2.}\quad For any $\sigma\in Aut{\cal L}$, we have that:
\vskip 0.2 cm (1)\quad $\sigma(L_0), \sigma(H_0)\in$ span$_{\ma F}\{
L_0, H_0, c \}$. \vskip 0.2 cm (2)\quad $ \sigma(G^+_\gamma)\in$
span$_{\ma F}\{G^+_\alpha\mid\alpha\in\Gamma \}$ or span$_{\ma
F}\{G^-_\alpha\mid\alpha\in\Gamma \}$ for all $\gamma\in\Gamma$.
\vskip 0.2 cm
 {\bf 
Proof.}\quad (1)\quad It is easy to see that $L_0, H_0, c,
G^\pm_{\alpha}$, for all $\alpha\in\Gamma$, are locally finite
elements of ${\cal L}$. And by lemma 3.1, we can prove (1).

(2) If there exist $\gamma_1, \gamma_2\in\Gamma$, satisfy that
$\sigma(G^+_{\gamma_1})\in$ span$_{\ma
F}\{G^+_\alpha\mid\alpha\in\Gamma \}$, $\sigma(G^+_{\gamma_2})\in$
span$_{\ma F}\{G^-_\alpha\mid\alpha\in\Gamma \}$, then we have that
\BC \sigma ([H_{\gamma_2-\gamma_1}, G^+_{\gamma_1}]) =
[\sigma(H_{\gamma_2-\gamma_1}), \sigma(G^+_{\gamma_1})]\in
\mbox{span}_{\ma F}\{G^+_\alpha \mid \alpha\in {\ma F}\}, \EC and
\BC \sigma ([H_{\gamma_2-\gamma_1}, G^+_{\gamma_1}]) =
\sigma(G^+_{\gamma_2})\in\mbox{span}_{\ma F}\{G^-_\alpha \mid
\alpha\in {\ma F}\},\EC it is a contradiction.
 \hfill$\Box$

\vskip 0.3 cm {\bf Proposition 3.3.}\quad For any $\sigma\in
Aut{\cal L}$, there exist $f\in$Hom$_{\ma Z}(\Gamma, {\ma F}^*),
\xi, \varepsilon\in\{\pm 1\}, a\in{\Gamma}, b\in{\ma F}^*$, such
that \BE\begin{array}{rll}\label{h1} \sigma(L_\alpha) & = &
\varepsilon f(\alpha)L_{\varepsilon\alpha}+f(\alpha)aH_{\varepsilon
a}
+\frac{1}{6}\varepsilon a^2 \delta_{\alpha,0}c, \\
[0.2 cm]\sigma(H_\alpha) & =
&\xi(f(\alpha)H_{\varepsilon\alpha}+\frac{a}{3}\varepsilon\delta_{\alpha,0}c),
\\ [0.2 cm] \sigma(G^+_{\alpha}) & = & f(\alpha)b(\delta_{\xi,1}
G^+_{\varepsilon(\alpha+a)}+\delta_{\xi,-1}G^-_{\varepsilon(\alpha-a)}),\\
[0.2 cm]\sigma(G^-_{\alpha}) & = & \varepsilon
f(\alpha)b^{-1}(\delta_{\xi,-1}
G^+_{\varepsilon(\alpha+a)}+\delta_{\xi,1}G^-_{\varepsilon(\alpha-a)}),\\
[0.2 cm] \sigma(c) & = & \varepsilon c,
\end{array}\EE where $\delta_{m,n}$ is also the Kronecker notation.

{\bf Proof.}\quad By lemma 3.1 and 3.2, we can assume that \BC
\sigma(L_0)=a_0L_0+b_0H_0+c_0c, ~~ \sigma(H_0)=d_0L_0+e_0H_0+f_0c,
\EC and \BC \sigma(L_\gamma) & = &
\sum_{\alpha\in\Gamma}a_{\gamma,\alpha}L_\alpha+\sum_{\alpha\in\Gamma}b_{\gamma,\alpha}H_\alpha+c_\gamma
c, \\  \sigma(H_\gamma) & = &
\sum_{\alpha\in\Gamma}d_{\gamma,\alpha}L_\alpha+\sum_{\alpha\in\Gamma}e_{\gamma,\alpha}H_\alpha+f_\gamma
c, \\   \sigma(c) & = & m_0c,\EC where $a_0, b_0, c_0, d_0, e_0,
f_0, m_0, a_{\gamma,\alpha}, b_{\gamma,\alpha}, c_{\gamma},
d_{\gamma,\alpha}, e_{\gamma,\alpha}, f_{\gamma} \in {\ma F}$.

{\bf Claim.} $a_0\not=0$.

If $a_0=0$, then $\sigma(L_0)=b_0H_0+c_0c$. Set $x\in{\cal
L}_\alpha$, $\alpha\not=0$, we have that
$-\alpha\sigma(x)=[\sigma(L_0), \sigma(x)]=[ b_0H_0+c_0c,
\sigma(x)]=0$, it is impossible. Therefore $a_0\not=0.$

 Acting $\sigma$ on $[L_0, L_\gamma]=-\gamma L_\gamma$,
we have that \BC -a_0\alpha
\sum_{\alpha\in\Gamma}(a_{\gamma,\alpha}L_\alpha+b_{\gamma,\alpha}H_\alpha)
=-\gamma\sum_{\alpha\in\Gamma}(a_{\gamma,\alpha}L_\alpha+b_{\gamma,\alpha}H_\alpha+c_\gamma
c), \EC i.e.,\BE\label{a2} \left\{\begin{array}{l} (\gamma-a_0\alpha
)a_{\gamma,\alpha}=0, \\ [0.2 cm] (\gamma-a_0\alpha
)b_{\gamma,\alpha}=0,
\\ [0.2 cm] c_{\gamma}=0. \end{array}\right.\EE
 If there exists $\gamma_0\in\Gamma$, such that
$\displaystyle{\frac{\gamma_0}{a_0}}\not\in\Gamma$,
 i.e., $\gamma_0-a_0\alpha\not=0$ for all
$\alpha\in\Gamma$, then $a_{\gamma_0, \alpha}=b_{\gamma_0,
\alpha}=0$ for all $\alpha\in\Gamma$, that is to say,
$\sigma(L_{\gamma_0})=0$, it is a contradiction since $\sigma$ is an
automorphism. Therefore,
$\displaystyle{\frac{\gamma}{a_0}}\in\Gamma$ for all
$\gamma\in\Gamma$, then $a_0=\pm 1$, where 1 is the unit of ${\ma
F}$.
\\ [0.3 cm] {\bf Case 1.}\quad $a_0=1$. \vskip 0.2 cm By (\ref{a2}),
we have that $\sigma(L_\gamma)=a_\gamma L_\gamma+b_\gamma
H_\gamma+c_0\delta_{\gamma,0}c$. And by the same argument, we also
have that $$\sigma(H_\gamma)=d_\gamma L_\gamma+e_\gamma
H_\gamma+f_0\delta_{\gamma,0}c.$$
 Applying $\sigma$ on $[H_\alpha,
H_\beta]=\frac{\alpha}{3}c\delta_{\alpha+\beta,0}$, we can obtain
that
$$
 \frac{\alpha}{3}m_0c\delta_{\alpha+\beta,0}=(\alpha-\beta)d_\alpha d_\beta L_{\alpha+\beta}
 +(\alpha d_\beta e_\alpha-\beta d_\alpha
 e_\beta)H_{\alpha+\beta}+(\frac{\alpha^3-\alpha}{12}d_{\alpha}d_{\beta}
 +\frac{\alpha}{3}e_\alpha e_\beta)c\delta_{\alpha+\beta,0},  $$
i.e.,$$ \left\{\begin{array}{l}(\alpha-\beta)d_\alpha d_\beta=0, \\
[0.2 cm] \alpha d_\beta e_\alpha-\beta d_\alpha e_\beta=0, \\ [0.2
cm] \frac{\alpha}{3}m_0=
\frac{\alpha^3-\alpha}{12}d_{\alpha}d_{-\alpha}
 +\frac{\alpha}{3}e_\alpha e_{-\alpha}. \end{array}\right.$$
If there exists $\alpha_0\in\Gamma$, such that $d_{\alpha_0}\not=0$,
then we will deduce that $d_\beta=e_\beta=0$ for all
$\beta\in\Gamma\setminus\{\alpha_0, 0\}$, contradiction with
$\sigma\in Aut{\cal L}$, hence $d_\alpha=0$ for all
$\alpha\in\Gamma$, and $\alpha m_0=\alpha e_\alpha e_{-\alpha}$,
i.e.,
$$d(H_\alpha)=e_\alpha H_\alpha+f_0\delta_{\alpha,0}c.$$
By lemma 3.2, we can discuss the action of $\sigma$ on ${\cal
L}_{\overline{1}}$ in two cases.
\\ [0.2 cm] {\bf Subcase 1.1.}\quad
$\sigma(G^+_\gamma)=\sum_{\alpha\in\gamma}g_{\gamma,\alpha}G^+_\alpha$.

Since $\sigma$ is an automorphism, by lemma 3.2, we have that
$\sigma(G^-_\gamma)=\sum_{\alpha\in\Gamma}h_{\gamma,\alpha}G^-_\alpha$.
By $\sigma[L_0, G^+_\gamma]=[\sigma(L_0), \sigma(G^+_\gamma)]$, we
get that
$$(\gamma+b_0-\alpha)g_{\gamma,\alpha}=0, $$ therefore,
$b_0\in\Gamma$, and $$\sigma(G^+_\gamma)=g_\gamma
G^+_{\gamma+b_0}.$$ Similarly, $$ \sigma(G^-_\gamma)=h_\gamma
G^-_{\gamma-b_0}.$$

  Applying
$\sigma$ to $[ H_\alpha, G^\pm_\beta ]$, we have that $$ \left\{
\begin{array}{lll} g_{\beta+\gamma}=e_\beta g_\gamma, \\ [0.2 cm]
h_{\beta+\gamma}=e_\beta h_\gamma,
\end{array}\right. ~~~ \mbox{for all} ~~ \beta, \gamma\in{\ma F}.
$$ Then we can deduce that
\BE\label{a4}  g_\alpha=e_\alpha g_0,  ~~~~~
 h_\alpha=e_\alpha h_0, ~~~~~ e_0=1. \EE
 Applying $\sigma$ to $[ G^-_\alpha,
G^+_\beta ]$, we have that \BE \label{a6} a_{\alpha+\beta} & = &
g_\alpha h_\beta, \\ [0.2 cm] \label{a8}
2b_{\alpha+\beta}-(\alpha-\beta)e_{\alpha+\beta} & = &
-(\alpha-\beta-2b_0)g_\alpha h_\beta, \\ [0.2 cm] \label{a10}
\frac{1}{3}((\alpha^2-b_0)^2-\frac{1}{4})g_\alpha h_{-\alpha} & = &
\frac{1}{3}(\alpha^2-\frac{1}{4})m_0-2\alpha f_0+2c_0. \EE By
(\ref{a4}) and (\ref{a6}), $$g_0h_0=a_0=1, $$ and $$
a_{\alpha+\beta}=e_\alpha g_0\cdot e_\beta h_0=e_\alpha e_\beta,
$$ then $a_\alpha=e_\alpha$ for all $\alpha\in{\Gamma}$.
Furthermore, following (\ref{a8}), we obtain that $$ b_\alpha=b_0
e_\alpha=b_0a_\alpha.
$$ By (\ref{a10}), let $\beta=-\alpha$, one can get that
$$ m_0=1, ~~~~ f_0=\frac{1}{3}b_0, ~~~~ c_0=\frac{1}{6}b_0^2. $$
 Therefore, \BE\begin{array}{llllll}\label{h3} \sigma(L_\alpha) & =
& a_\alpha L_{\alpha}+a_\alpha b_0H_{\alpha} +\frac{1}{6}b_0^2 c
\delta_{\alpha,0}, &  \sigma(H_\alpha) & = & a_\alpha
H_{\alpha}+\frac{1}{3}b_0 c\delta_{\alpha,0},
\\ [0.2 cm] \sigma(G^+_{\alpha}) & = & a_\alpha g_0
G^+_{\alpha+b_0},&  \sigma(G^-_{\alpha}) & = & a_\alpha g_0^{-1}
G^-_{\alpha-b_0},\\
[0.2 cm] \sigma(c) & = &  c,
\end{array}\EE where $a_\alpha\in{\ma F^*}$ satisfies that $a_{\alpha+\beta}=a_\alpha a_\beta$, $b_0\in \Gamma, g_0\in{\ma
F^*}$.
\\ [0.2 cm] {\bf Subcase 1.2.}\quad
$\sigma(G^+_\gamma)=\sum_{\alpha\in\gamma}u_{\gamma,\alpha}G^+_\alpha$.
\vskip 0.2 cm  Under the same argument, we can obtain that
\BE\begin{array}{llllll}\label{h5} \sigma(L_\alpha) & = & a_\alpha
L_{\alpha}+a_\alpha b_0H_{\alpha} +\frac{1}{6}b_0^2 c
\delta_{\alpha,0}, & \sigma(H_\alpha) & = & -a_\alpha
H_{\alpha}-\frac{1}{3}b_0 c\delta_{\alpha,0},
\\ [0.2 cm] \sigma(G^+_{\alpha}) & = & a_\alpha u_0
G^-_{\alpha-b_0},&  \sigma(G^-_{\alpha}) & = & a_\alpha u_0^{-1}
G^+_{\alpha+b_0},\\
[0.2 cm] \sigma(c) & = &  c,
\end{array}\EE where $a_\alpha\in{\ma F^*}$ satisfies that $a_{\alpha+\beta}=a_\alpha a_\beta$, $b_0\in \Gamma, u_0\in{\ma
F^*}$. \\ [0.2 cm]{\bf Case 2.}\quad $a_0=-1$.

In this case, we have that $\sigma(L_0)=-L_0+b_0H_0+c_0c$. By the
similar method, we can obtain the two following subcases:
\\
[0.2 cm] {\bf Subcase 2.1.} \BE\begin{array}{llllll}\label{h7}
\sigma(L_\alpha) & = & -a_\alpha L_{-\alpha}+a_\alpha b_0H_{-\alpha}
-\frac{1}{6}b_0^2 c\delta_{\alpha,0}, &  \sigma(H_\alpha) & = &
a_\alpha H_{-\alpha}-\frac{1}{3}b_0 c\delta_{\alpha,0},
\\ [0.2 cm] \sigma(G^+_{\alpha}) & = & a_\alpha g_0
G^+_{-\alpha-b_0},&  \sigma(G^-_{\alpha}) & = & -a_\alpha g_0^{-1}
G^-_{-\alpha+b_0},\\
[0.2 cm] \sigma(c) & = &  -c.
\end{array}\EE
 {\bf Subcase 2.2.}\BE\begin{array}{llllll}\label{h9} \sigma(L_\alpha) & = & -a_\alpha
L_{-\alpha}+a_\alpha b_0H_{-\alpha} -\frac{1}{6}b_0^2 c
\delta_{\alpha,0}, & \sigma(H_\alpha) & = & -a_\alpha
H_{\alpha}+\frac{1}{3}b_0 c\delta_{\alpha,0},
\\ [0.2 cm] \sigma(G^+_{\alpha}) & = & a_\alpha u_0
G^-_{-\alpha+b_0},&  \sigma(G^-_{\alpha}) & = & -a_\alpha u_0^{-1}
G^+_{-\alpha-b_0},\\
[0.2 cm] \sigma(c) & = & -c,
\end{array}\EE where $a_\alpha\in{\ma F^*}$ satisfies that $a_{\alpha+\beta}=a_\alpha a_\beta$, $b_0\in \Gamma, g_0, u_0\in{\ma
F^*}$. \hfill $\Box$ \vskip 0.3 cm
 By proposition 3.3, we can obtain
that
$$ Aut{\cal L}=\{\sigma(f, \xi, \varepsilon, a, b)\mid f(\alpha)\in\mbox{Hom}_{\ma Z}(\Gamma, {\ma F^*}),
 \xi, \varepsilon\in\{\pm 1\},
 a\in \Gamma, b\in{\ma
F^*}\}, $$ where $\sigma(f, \xi, \varepsilon, a, b)\in Aut{\cal L}$
 satisfies the following relations:
\BE\begin{array}{lll}\label{h11} \sigma(f, \xi, \varepsilon, a,
b)(L_\alpha) & = & \varepsilon f(\alpha)
L_{\varepsilon\alpha}+f(\alpha) a H_{\varepsilon\alpha}
+\frac{1}{6}\varepsilon a^2 c\delta_{\alpha,0}, \\ [0.2 cm]
\sigma(f, \xi, \varepsilon, a, b)(H_\alpha) & = & \xi (f(\alpha)
H_{\varepsilon\alpha}+\frac{1}{3}\varepsilon a c\delta_{\alpha,0}),
\\ [0.2 cm] \sigma(f, \xi, \varepsilon, a, b)(G^+_{\alpha}) & = &
f(\alpha) b (\delta_{\xi,1}
G^+_{\varepsilon({\alpha+a})}+\delta_{\xi,-1}
G^-_{\varepsilon({\alpha-a})}),\\ [0.2 cm] \sigma(f, \xi,
\varepsilon, a, b)(G^-_{\alpha}) & = & \varepsilon f(\alpha) b^{-1}
(\delta_{\xi,-1} G^+_{\varepsilon({\alpha+a})}+\delta_{\xi,1}
G^-_{\varepsilon({\alpha-a})}),\\
[0.2 cm] \sigma(f, \xi, \varepsilon, a, b)(c) & = & \varepsilon c,
\end{array}\EE where  $f(\alpha)\in$Hom$_{\ma Z}(\Gamma, {\ma F^*})$, i.e.,
$f(\alpha+\beta)=f(\alpha)f(\beta)$ for all $\alpha,
\beta\in\Gamma$, and $\xi, \varepsilon\in\{\pm 1\}$,
 $a\in \Gamma, b\in{\ma
F^*}$. Obviously, $\sigma(f, 1, 1, a, b)$, $\sigma(f, -1, 1, a, b)$,
$\sigma(f, 1, -1, a, b)$) and $\sigma(f, -1, -1, a, b)$ are the
automorphisms which are defined in (\ref{h3})-(\ref{h9}). And
$\sigma(f_1, \xi_1, \varepsilon_1, a_1, b_1)=\sigma(f_2, \xi_2,
\varepsilon_2, a_2, b_2)$ if and only if $f_1=f_2$, $\xi_1=\xi_2,
\varepsilon_1=\varepsilon_2$, $a_1=a_2, b_1=b_2$. \BC & &
\sigma(f_1, \xi_1, \varepsilon_1, a_1, b_1)\sigma(f_2, \xi_2,
\varepsilon_2, a_2, b_2) \\ [0.2 cm] & = & \sigma(\varphi,
\xi_1\xi_2, \varepsilon_1\varepsilon_2,
\xi_1\xi_2\varepsilon_1a_1+\xi_2a_2,f_1(\xi_2\varepsilon_2a_2)b_1^{\xi_2}b_2),\EC
where $\varphi$ satisfies that
$\varphi(\alpha)=f_1(\varepsilon_2\alpha)f_2(\alpha)$ for any
$\alpha\in{\ma F}$. Obviously, $Aut{\cal L}$ is not a ablian group.
For any $f\in$Hom$_{\ma Z}(\Gamma, {\ma F})$, we define $f^{-1}$ by
$f^{-1}(\alpha)=(f(\alpha))^{-1}$ for any $\alpha\in\Gamma$. It is
easy to see that $f^{-1}$ is also a homomorphism from $\Gamma$ to
${\ma F}^*$. Then for any $\sigma(f, \xi, \varepsilon, a, b)\in
Aut{\cal L}$, we can get that $$\sigma^{-1}(f, \xi, \varepsilon, a,
b)=\sigma(f^{-\varepsilon}, \xi, \varepsilon, -\xi\varepsilon a,
sgn(\xi+\varepsilon)f(a)b^{-1}),$$ where
$sgn(\eta)=\left\{\begin{array}{ll} 1, & \eta\geqslant 0 \\  -1, &
\eta<0
\end{array}\right.$, and $f^1=f$.

Set $\tau=<\sigma(f, 1, 1, a, b)>$, it is not difficult to prove
that $\tau$ is a normal subgroup of $Aut{\cal L}$. Then we can
obtain the following main theorem:
 \vskip 0.3 cm {\bf Theorem 3.4.}\quad (1) $Aut{\cal L}$ is isomorphic to
  $\mbox{Hom}_{\ma Z}(\Gamma, {\ma F^*})\times
 {\ma Z'}_2\times{\ma Z'}_2\times\Gamma\times{\ma F^*}$, where ${\ma Z'}_2=\{\pm 1\}$, the group multiplication is given by:
$$ (f_1, \xi_1, \varepsilon_1, a_1, b_1)\cdot (f_2, \xi_2, \varepsilon_2, a_2, b_2)=(\varphi, \xi_1\xi_2,
\varepsilon_1\varepsilon_2,
\xi_1\xi_2\varepsilon_1a_1+\xi_2a_2,f_1(\xi_2\varepsilon_2a_2)b_1^{\xi_2}b_2),$$
where $\varphi$ satisfies that
$\varphi(\alpha)=f_1(\varepsilon_2\alpha)f_2(\alpha)$ for any
$\alpha\in{\ma F}$.
   \\ [0.2 cm]
 (2)  $Aut{\cal
L}/\tau$ is a Klein group.

    \end{document}